\documentclass[12pt]{article}
\usepackage{pstricks}
\usepackage{graphicx}
\usepackage{amsmath}
\usepackage{amssymb}
\usepackage{amsthm}
\makeatletter
\renewcommand{\@makecaption}[2]{%
\vspace{\abovecaptionskip}%
\sbox{\@tempboxa}{#1. #2}
#1. #2\par
\vspace{\belowcaptionskip}} 
\makeatother
\begin{document}

\begin{center}

{\LARGE \bf Weak forms of shadowing in topological dynamics}

\bigskip

{\noindent\large Danila Cherkashin${}^{a,b,c}$, Sergey Kryzhevich${}^{a,d}$\footnote{Corresponding author.\\\hspace*{5.5mm}Email address:
kryzhevicz@gmail.com}}

\medskip

{\noindent\small
${}^a$ Faculty of Mathematics and Mechanics, Saint-Petersburg State University,\\
 28, Universitetskiy pr., Peterhof, Saint-Petersburg, 198503, Russia;\\
${}^b$ Faculty of Innovations and High Technology, Moscow Institute of Physics and Technology, \\ 9, Institutskiy per., Dolgoprudny,
Moscow Region, 141700,
Russia;\\
${}^c$ Section of Mathematics, Geneva University, \\ 2-4 rue du Li\'evre,
Case postale 64
1211 Geneva 4, Switzerland;\\
${}^d$ School of Natural Sciences and Mathematics, The University of Texas at Dallas,\\ 800 W. Campbell Road, Richardson, Texas, 75080-3021, USA
}
\end{center}

\begin{abstract} An approach to find a weak form of shadowing is developed.  We consider continuous maps of compact metric spaces. It is proved that every pseudotrajectory with sufficiently small errors contains at least one subsequence that can be shadowed by a subsequence of an exact trajectory with same indices. Later, we study homeomorphisms such that any pseudotrajectory can be shadowed by a finite number of exact orbits. We call this property multishadowing. Criteria for existence of $\varepsilon$-- networks whose iterations are $\varepsilon$ -- networks are given. Relations between multishadowing and some ergodic and topological properties of dynamical systems are discussed. Various applications of obtained results are given.
\end{abstract}

\textbf{Keywords:}
topological dynamics, minimal points, invariant measure, shadowing, chain recurrence, $\varepsilon$ -- networks, syndetic sets.

\section{Introduction}
Shadowing is a very important property of dynamical systems, closely related to problems of structural stability and modeling. For review on general Shadowing Theory we refer to [21,26--28].

Though the most evident application of shadowing is related to numerical methods, first results involving the concept of pseudotrajectories were obtained by Anosov [2], Bowen [8] and Conley [10] as a tool to study qualitative properties of dynamical systems.

In a nutshell, shadowing is existence of an exact trajectory pointwise near a given pseudotrajectory i.e. a trajectory with errors. This property is closely related to structural stability. Indeed, it is well-known that structural stability implies shadowing [34,38]. Such shadowing is Lipschitz [29].

Sakai [36] demonstrated that the $C^1$ -- interior of the set of all diffeomorphisms with shadowing coincides with the set of all structurally stable diffeomorphisms. Osipov, Pilyugin and Tikhomirov [25,29] demonstrated that the so-called Lipschitz periodic shadowing property is equivalent to $\Omega$ -- stability, see also [27]. Moreover, the corresponding set of dynamical systems coincides with the interior of the set of systems with periodic shadowing property and with the set of systems with orbital limit shadowing property.

Pilyugin and Tikhomirov [33] demonstrated that Lipschitz shadowing is equivalent to structural stability.

Shadowing is not $C^1$ generic. Bonatti, Diaz and Turcat [7] demonstrated that there is a $C^1$ -- open set of diffeomorphisms of the 3--torus where none of diffeomorphisms satisfies shadowing property. Yuan and Yorke [41] proved a similar result for $C^r$ -- diffeomorphisms ($r>1$).

Surprisingly, shadowing is generic in the $C^0$ topology of homeomorphisms of a smooth manifold. This was proved by Pliyugin and Plamenevskaya [30]. Similar results were obtained  for continuous mappings of manifolds [19] and for continuous maps of Cantor set [4].

This fact inspires studying shadowing by means of topological dynamics. This approach gave many important results mostly obtained in last two decades.

Mai and Ye [22] demonstrated that odometers have shadowing. This is the only example of such type infinite minimal systems. Of course, there are many non-minimal infinite systems with shadowing e.g. Bernoulli shift.

On the other hand, Moothathu [23] proved that minimal points are dense for every non-wandering system with shadowing. Moothathu and Oprocha [24] demonstrated that non-wandering systems with shadowing have a dense set of regularly recurrent points. 

Dastjerdi and Hosseini [11,12] studied "almost identical"\ mappings. They proved that if a chain transitive dynamical system has an equicontinuity point then it is a distal, equicontinuous and minimal homeomorphism (see also [15,16]). Thus any transitive system with shadowing is either sensitive or equicontinuous.

Another version of shadowing (the so-called average shadowing) was introduced by Blank [5]. The so-called ergodic shadowing was studied in [13]. Some other kinds of shadowing ($\underline d$--shadowing, weak shadowing, etc.) were discussed in [12, 35] and [37], see also references therein.

However, the problem of shadowing in non-smooth dynamical systems is very far from being resolved. Theoretical results in this area may be applied for modeling non-smooth dynamics like vibro-impact systems, systems with dry friction, biological problems etc.  

The main objective of this paper is the following. We demonstrate that for a very general dynamical system, any numerical method, even an inappropriate one, can give some useful information on asymptotical behavior of solutions. First of all, it can be used to find an invariant measure (Theorem 3.1). If we take a random point of a pseudotrajectory, obtained by this "incorrect"\ numerical method, there is a positive probability to find a minimal point in a neighborhood of the selected point (Theorem 3.1, Corollary 4.5). In some generic assumptions (see Theorem 3.3) this probability is equal to 1.

First of all, we show that for any dynamical system and any pseudotrajectory there is a subsequence that can be traced by a subsequence of a precise trajectory with same indices. This is the first key result of our paper.

Then it is natural to ask, if any pseudotrajectory can be traced by a finite number of trajectories. This is the so-called multishadowing (Definition 2.18). We demonstrate that this property is $C^1$ -- generic, we describe it in terms of topological characteristics of dynamical systems e.g. minimal points, invariant measures, etc.  We study a generalization of equicontinuous systems i.e. systems with almost invariant $\varepsilon$ -- networks. An $\varepsilon$ -- network is called almost invariant if all its iterations are $\varepsilon$ -- networks.

The second central statement of our research is Theorem 3.3. We prove that for a nonwandering system multishadowing is equivalent to existence of almost invariant $\varepsilon$ -- networks for any $\varepsilon>0$. Moreover, both these properties are equivalent to the so-called Bronstein condition [9] i.e. density of minimal points in the set of nonwandering points (Definition 2.8).

The paper is organized as follows. First of all, we recall the terminology, related to Shadowing Theory and Topological Dynamics (Section 2). In Section 3 we list principal results of the paper.

We improve the main result of [20] in Section 4. It is proved that for any continuous mapping of a compact metric space into itself and for any one-side pseudotrajectory $x_k$ there exists a sequence $k_n$ and a precise trajectory 
$\{y_k=T^k(y_0)\}$ such that points $x_{k_n}$ and $y_{k_n}$ are uniformly close. The density of $\{k_n\}$ in $\mathbb N$ is positive (Theorem 3.1).

In Sections 5 and 6 we study nonwandering systems. We prove that multishadowing is equivalent to Bronstein condition. In Section 7 we prove that multishadowing is equivalent to existence of almost invariant $\varepsilon$ -- networks for all $\varepsilon>0$. Moreover, for nonwandering homeomorphisms, multishadowing implies existence of an invariant measure, supported on all the phase space (Section 8). 

In Section 9 we prove that if every chain recurrent point is nonwandering and Bronstein condition holds on the nonwandering set, the considered system satisfies multishadowing property. The converse statement is proved in Section 10.

In Section 11 we study networks that are almost invariant almost everywhere with respect to an invariant measure. 

In Section 12 we demonstrate that multishadowing is $C^0$ and $C^1$ -- generic.

In Sections 13 we discuss possible applications of the main results of the paper.

Conclusion is given in Section 14.

\section{Definitions}

Recall some standard definitions from Topological Dynamics. Consider a compact metric space $X$ endowed with the metric $\rho$.  Let a map $T:X\to X$ be continuous.

\noindent\textbf{Definition 2.1.} Let $d>0$. A sequence $\{x_k\}_{k\in {\mathbb N}}$ is a \emph{$d$ -- pseudotrajectory} if
$$\rho(x_{k+1}, T(x_k))\le d$$ for all $k \in {\mathbb N}$.

\noindent\textbf{Definition 2.2.} We say that the mapping $T$ satisfies \emph{shadowing property} if for any $\varepsilon>0$ there is a $d>0$ such that for any $d$ -- pseudotrajectory $\{x_k\}$ there exists an exact trajectory $\{y_k=T^k(y_0), k\in {\mathbb N}\}$ such that
$$\rho(x_k,y_k)<\varepsilon \eqno (2.1)$$
for all $k\in {\mathbb N}$.

If $T:X\to X$ is a homeomorphism, we may consider "two-sided"\ pseudotrajectories $\{x_k\}_{k\in {\mathbb Z}}$ and study "two-sided shadowing", defined similarly to Definition 2.2. Abusing notations, we say "pseudotrajectory"\ and "shadowing"\ in both cases. If it is necessary we add words "one-sided"\ or "two-sided"\ in order to underline which kind of dynamical systems we deal with.

\noindent\textbf{Definition 2.3.} A point $x\in X$ is \emph{wandering} if there exists a neighborhood $U\ni x$ such that $T^k(U)\bigcap U=\emptyset$ for all $k\in {\mathbb N}$.

\noindent\textbf{Definition 2.4.} \emph{Non-wandering} points form the \emph{non-wandering set} $\Omega(X,T)$. Let $\mathrm{NW}$ be the class of non-wandering systems ($X=\Omega (X,T)$).

\noindent\textbf{Definition 2.5.} A point $y\in X$ is an $\omega$ -- limit point for $x\in X$ i.e. $y\in \omega(x)$ if there exists a sequence $n_k\to +\infty$ such that $T^{n_k} (x)\to y$. Let $\omega(X,T)$ be the closure of all $\omega$ -- limit points for all points of $X$.

Recall some classic notations. Define the positive semiorbit of a point $x$ by formula $O^+(x)=\{T^k(x):k\ge 0\}$. For homeomorphisms, we consider orbits: $O(x)=\{T^k(x):k\in {\mathbb Z}\}$.

\noindent\textbf{Definition 2.6.} The dynamical system $(X,T)$ is called \emph{minimal}, if $\overline{O^+(x)} = X$ for every $x \in X$. 

\noindent\textbf{Definition 2.7.} A point $y\in X$ is called \emph{minimal} (or almost periodic) for dynamical system $(X,T)$, if the subsystem $(\overline{O^+(y)},T)$ is minimal. Let $\mathrm{M}(X,T)$ be the set of all minimal points of $(X,T)$.

\noindent\textbf{Definition 2.8.} If the set of minimal points is dense in $X$ we say that $(X,T)$ satisfies the \emph{Bronstein condition.}

Let us also recall a definition from Combinatorics and Number Theory.

\noindent\textbf{Definition 2.9.} A subset $S\subset {\mathbb N}$ is called \emph{syndetic} if there exists $n=n(S)\in {\mathbb N}$ such that for any $m\in {\mathbb N}$ the intersection $S\bigcap [m,m+n]$ is non-empty. We also use notion of $n$ -- syndetic set if we need to specify the value $n$. 

We recall a well-known fact from the theory of minimal sets [14,18].

\noindent\textbf{Lemma 2.10.} \emph{Let $T:X\to X$ be a continuous map. System $(X, T)$ is minimal if and only if the set
$$
N(x, U)=\{m\in {\mathbb N}: T^m(x)\in U\} \eqno (2.2)
$$
is syndetic for every $x \in X$ and nonempty open set $x\in U \subset X$.}

Starting form here we assume up to the end of the section that $T:X\to X$ is a homeomorphism.

\noindent\textbf{Definition 2.11.} We say that a point $z\in X$ is an $\alpha$ -- limit point for a point $x\in X$ if there exists an integer sequence $n_k\to \infty$ such that $T^{-n_k} (x)\to z$. Let $\alpha(X,T)$ be the closure of all $\alpha$ -- limit points for all points of $X$.

\noindent\textbf{Definition 2.12.} A point $x\in X$ is \emph{recurrent} if $x\in \alpha(x)\bigcap \omega(x)$. Let $\mathrm{R}(X,T)$ be the set of all recurrent points of system $(X,T)$.

\noindent\textbf{Definition 2.13.} The \emph{chain recurrent set} $\mathrm{CR} (X, T)$ is the set of points $x \in X$ such that for any $d>0$ there exists a finite $d$ -- pseudotrajectory $x=x_1, \dots, x_k=x$, $k>1$.

We recall a well-known result from Topological Dynamics.

\noindent\textbf{Lemma 2.14.} \emph{Let $T:X\to X$ be a homeomorphism, $\mathrm{P}(X,T)$ be the set of all periodic points of $T$. Then
\begin{enumerate}
\item sets $\Omega(X,T)$ and $\mathrm{CR}(X,T)$ are closed;
\item $\mathrm{P}(X,T)\subset \mathrm{M}(X,T)\subset \mathrm{R}(X,T) \subset \alpha(X,T)\bigcup \omega(X,T) \subset \Omega (X,T) \subset \mathrm{CR}(X,T)$;
\item \textrm{[17, Proposition 4.1.18]} if $\mu$ is a Borel probability invariant measure for $(X,T)$ then $\mathrm{supp}\, \mu\subset \overline{\mathrm{R}(X,T)}$.
\end{enumerate}}

Here we recall that the support $\mathrm{supp}\, \mu$ of a Borel measure $\mu$ is the intersection of all closed subsets $Y\subset X$ such that $\mu(Y)=1$. 

\noindent\textbf{Definition 2.15.} A subset $Y\subset X$ is an $\varepsilon$ -- network in $X$ if for any $x\in X$ there exists a $y\in Y$ such that $\rho(x,y)\le \varepsilon$. 

\noindent\textbf{Definition 2.16.} An $\varepsilon$-network $Y$ is \emph{almost invariant} if for every $n \in \mathbb{Z}$ the set $T^n(Y)$ is an $\varepsilon$ -- network. (Fig.\, 1).

\begin{figure}[ht!]
\begin{center}
\includegraphics*[width=4in]{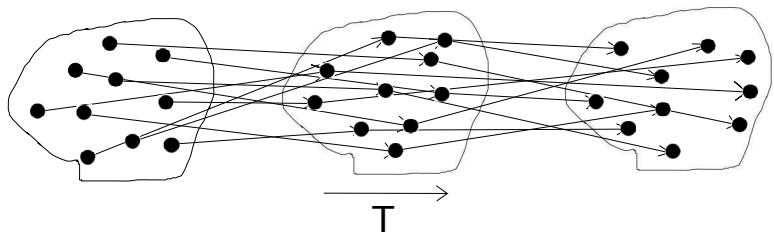}
\end{center}
\caption{Almost invariant networks.}
\end{figure}

Denote by $\mathrm Q$ the class of systems $(X,T)$ that have finite almost invariant $\varepsilon$-networks for every $\varepsilon>0$.

\noindent\textbf{Lemma 2.17}. $\mathrm{Q}\subset \mathrm{NW}$. 

\noindent\textbf{Proof.}
If $(X,T)\in \mathrm{Q}$ any neighborhood of any point of $X$ contains an $\omega$ -- limit point, corresponding to a limit point of one of points of an almost invariant network. $\square$

\noindent\textbf{Definition 2.18} (Fig.\,2).
Let $\mathrm{W}$ be the class of dynamical systems $(X,T)$ such that for any $\varepsilon>0$ there exists a $d>0$ as follows: for any $d$ -- pseudotrajectory $\{x_k\}$ there exist points $y^1,\ldots, y^N$ ($N$ may depend on $\{x_k\}$ and $\varepsilon$) such that $x_k$ is $\varepsilon$ close to one of points $T^k(y^i)$ for all $k\in {\mathbb N}$.
Then the system $(X,T)$  is said to satisfy the \emph{multishadowing property}. The corresponding maximal number of shadowing trajectories $N(\varepsilon)$ is called \emph{multishadowing parameter}.

\begin{figure}
\begin{center}
\includegraphics*[width=4in]{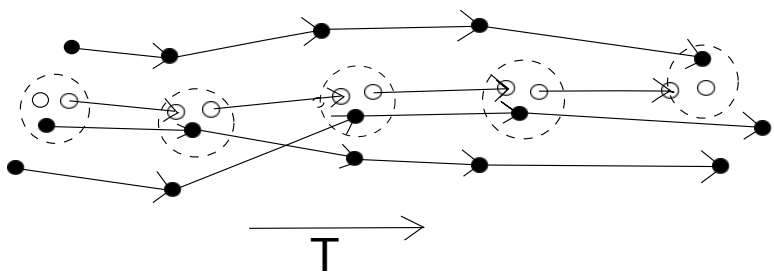}
\end{center}
\caption{Multishadowing.}
\end{figure}

Of course, shadowing implies multishadowing. The converse statement is not true. For instance, 
$(X,\mathrm{id})\in \mathrm{W}$ for any compact metric space $X$. Another counterexample, one may keep in mind, is a discretization of the o.d.e. $\dot x=x^2-x^4$, defined on the segment $[-1,1]$ (Fig.\, 3).

\begin{figure}
\begin{center}
\includegraphics*[width=2.5in]{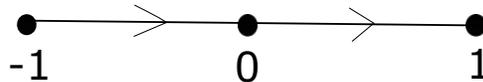}
\end{center}
\caption{A system with multishadowing and without shadowing.}
\end{figure}

There exist dynamical systems that do not belong to the class $\mathrm W$ (see Proposition 4.3 below).

\noindent\textbf{Definition 2.19.} We say that system $(X,T)$ is \emph{equicontinuous} if the family of maps $T^k:X\to X$, $k\in {\mathbb Z}$ is equicontinuous.

\noindent\textbf{Remark 2.20.} The class $Q$ is a natural generalization of equicontinuous systems. Evidently, all equicontinuous systems belong to $Q$. Meanwhile, the introduced class is much reacher, it includes some expansive systems e.g. dynamics on non-wandering sets for Axiom A diffeomorphisms.

\noindent\textbf{Definition 2.21.} Let $T$ be a homeomorphism of a compact metric space $X$, $\mu$ be a Borel probability invariant measure on $X$. We say that a finite set $A$ is an $\varepsilon$ -- network, \emph{almost invariant with respect to $\mu$} if $\mu(U_{\varepsilon}(T^n(A)))>1-\varepsilon$ for any $n\in {\mathbb Z}$.

Here $U_\varepsilon (\cdot)$ stands for $\varepsilon$ neighborhood of a set. 

\section{Main results}

\noindent\textbf{Theorem 3.1.} \emph{Let $T:X\to X$ be a continuous map of a compact metric space $X$. For any $\varepsilon>0$ there exists a $d>0$ such that for any one-side $d$ -- pseudotrajectory $\{x_k, k\ge 0\}$ there exists a subsequence $\{k_n\}$ and a point $y\in M(X,T)$ such that $\rho(x_{k_n}, T^{k_n}(y))<\varepsilon$.
The sequence $k_n$ may be taken so that
$$a:=\limsup\limits_{N\to \infty} \dfrac{\# K\bigcap [0,N]}{N}>0. \eqno (3.1)$$
Here $K=\{k_n,\, n\in {\mathbb N}\}$.}

If Eq. (3.1) is satisfied, we say that the set $K$ has \emph{positive density} in ${\mathbb Z}^+$.

This result is proved and discussed in Section 4. In fact we do not prove that for a given pseudotrajectory there is a trajectory that traces it. We just prove that both the pseudotrajectory and the "shadowing"\ trajectory return to a neighborhood of the same point along the same sequence of instants of time.

\noindent\textbf{Remark 3.2.} A result very similar to Theorem 3.1 was proved by one of co-authors in [20]. However, the statement of Theorem 3.1 is stronger. In [20] it was not proved that the sequence $\{k_n\}$ can be chosen so that (3.1) is satisfied.  In other words, we prove that the sequence $\{k_n\}$ does not grow too fast, that may be important for applications. In order to obtain inequality (3.1) we have to modify the proof (see Section 4).

Let $\mathrm{Br}$ be the class of all systems, corresponding to homeomorphisms of $X$ that satisfy Bronstein condition (see Definition 2.8). Recall that $\mathrm W$ is the class of dynamical systems with the muslishadowing property and 
$\mathrm Q$ is the class of systems that have almost invariant $\varepsilon$ -- networks for all $\varepsilon>0$. 

\noindent\textbf{Theorem 3.3.} \emph{
\begin{enumerate}
\item $\mathrm{Q} = \mathrm{Br} = \mathrm{W} \bigcap \mathrm{NW}$.
\item For any homeomorphism from the class $\mathrm Q$ there exists a probability invariant measure, supported on all $X$.
\item
$(X, T) \in\mathrm W$
if and only if 
$$\mathrm{CR} (X,T) =\overline{\mathrm{M}(X,T)}. \eqno (3.2)$$
\end{enumerate}}

\noindent\textbf{Remark 3.4.}
It is more convenient for us to deal with the following conditions, equivalent to $(3.2)$:
\begin{enumerate}
\item
$$\mathrm{CR} (X,T) =\Omega(X,T); \eqno (3.3)$$
\item Bronstein condition holds for system $(\Omega(X,T),T)$.
\end{enumerate}

We split the statement of Theorem 3.3 to several lemmas.

\noindent{\bf Lemma 3.5.} {\it Systems $(X,T)$ that satisfy Bronstein condition belong to the class $\mathrm Q$.}

\noindent{\bf Lemma 3.6.} {\it Let $K$ be a compact invariant set for system $(X,T)$. Assume that for any $\varepsilon>0$ there exists a finite set $A_\varepsilon\subset X$ such that $K\subset U_\varepsilon (T^k(A_\varepsilon))$ for any $k\in {\mathbb Z}$. Then $K\subset \overline{\mathrm{M}(X,T)}$.}

Observe that here we do not assume that $A_\varepsilon \subset K$. Taking $K=X$, we obtain $\mathrm{Q}\subset \mathrm{Br} \bigcap \mathrm{NW}$.

\noindent{\bf Lemma 3.7.} {\it $\mathrm{Q} \subset \mathrm{W}$; $(X,T)\in W$ implies $(\mathrm{CR}(X,T),T)\in \mathrm{Q}$.}

Particularly, $\mathrm{W}\bigcap \mathrm{NW}\subset \mathrm{Q}$. So, the first part of Theorem 3.3 follows from Lemmas 3.5 -- 3.7.

\noindent{\bf Lemma 3.8.} {\it If $(X,T)\in \mathrm Q$, there exists of Borel probability invariant measure, supported on all $X$.}

By virtue of [17,Theorems 4.1 and 7.1] existence of such an invariant measure implies that $X=\overline{\mathrm{R}(X,T)}$.

\noindent{\bf Lemma 3.9. } {\it Let (3.2) take place. Then system $(X,T)$ has multishadowing property.}

\noindent{\bf Lemma 3.10. }{\it $(X,T)\in \mathrm W$ implies Eq.\, $(3.3)$.}

Statements of Lemmas 3.9 and 3.10 imply the third item of Theorem 3.3.

Finally, we formulate an "ergodic"\ version of Lemmas 3.5 and 3.6.

\noindent\textbf{Theorem 3.11.} \emph{Let $T$ be a homeomorphism of a compact metric space $X$, $\mu$ be a Borel probability invariant measure on $X$. 
\begin{enumerate}
\item If for any $\delta>0$ there exists a finite $\delta$ -- network $A_\delta$, almost invariant with respect to $\mu$, then $\mathrm{supp}\, \mu\subset \overline{{\mathrm M}(X,T)}$.
\item If $\mathrm{supp}\,\mu \subset \overline{{\mathrm M}(X,T)}$ we can take $A_\varepsilon\subset \mathrm{supp}\,\mu$ for any $\varepsilon>0$. 
\end{enumerate}}

\noindent\textbf{Remark 3.12.} Density of minimal points for nonwandering systems with shadowing was proved by Moothathu [23,Theorem 1]. Theorems 3.5 and 3.6 demonstrate that so it is for nonwandering systems with multishadowing. 

\section{Partial shadowing. Proof of Theorem 3.1}

First, we prove an auxiliary statement. 

\noindent\textbf{Lemma 4.1.} \emph{For any positive sequence $\delta_m\to 0$ and for any sequence $\{p_k^m\}$ of $\delta_m$ -- pseudotrajectories there exists a point $\bar x\in \mathrm{M}(X,T)$ such that sets 
$$S_m=\{k:p_k^m\in B_{\varepsilon/2} (\bar x)\}$$ 
where $m$ is sufficiently big have positive densities in ${\mathbb Z}^+$.}

\noindent\textbf{Proof of Lemma 4.1.} We use some ideas of the proof of the Krylov-Bogolyubov Theorem [17, Theorem 4.1.1].

Fix corresponding sequences $\delta_m$ and $p_k^m$. 

Let $\Phi=\{\varphi_k:k\in {\mathbb N}\}$ be a countable sets of continuous functions on $X$, dense in $C^0(X\to {\mathbb R})$. 

Using diagonal sequence method, we obtain an integer sequence $s_j\to \infty$ such that for any function $\varphi\in \Phi$ there exists a limit
$$J_m(\varphi):=\lim_{j\to\infty} \dfrac{1}{s_j}\sum\limits_{i=0}^{s_j-1} \varphi(p^m_i). \eqno (4.1)$$
Moreover, we can take the diagonal sequence so that the sequence $\{s_j\}$ is the same for all $m$.

Let us demonstrate that functionals $J_m$ can be continuously extended to $C^0(X\to {\mathbb R})$.

Indeed, let $\psi \in C^0(X\to {\mathbb R})$ and $\varepsilon>0$. Take a function $\varphi\in \Phi$ so that $\|\psi-\varphi\|_{C^0}\le \varepsilon$. Then, for any $j\in {\mathbb N}$ we have
$$\left| \dfrac{1}{s_j}\sum\limits_{i=0}^{s_j-1} (\varphi(p^m_i)-\psi(p^m_i))\right|\le \varepsilon.$$

This demonstrates that the value $J_m(\psi)$ is correctly defined by the formula, similar to (4.1). Moreover, $|J_m (\psi)|\le \|\psi\|_{C^0}$.

Evidently $\|J_m\|\le 1$ for all $m$. So, all functionals $J_m:C(X \to {\mathbb R})$ are linear, continuous and positive. By virtue of Riesz Representation Theorem, they uniquely define probability measures $\mu_m$ on $X$ according to the formula
$$J_m(\varphi)=\int_X \varphi \, d\mu_m. \eqno (4.2)$$

By virtue of Banach-Alaoglu Theorem, the set of all Borel probability measures is compact in the $*$-weak topology.
Without loss of generality, we can suppose that the considered sequence $*$ -- weakly converges to a Borel probability measure $\mu_*$. Let us demonstrate that $\mu_*$ is an invariant measure. Fix a $\varphi\in C(X)$, then
$$\int_X ((\varphi\circ T)-\varphi )\, d\mu_*=\lim_{m\to\infty} (J_m(\varphi\circ T)-J_m (\varphi))=$$
$$\lim_{m\to\infty}\lim_{j\to\infty} \dfrac{1}{s_j}\left(\sum\limits_{i=1}^{s_j-1} (\varphi(T(p^m_{i-1}))-\varphi(p_i^m))+\varphi(T(p^m_{s_j-1}))-\varphi(p^m_0)\right)=0. \eqno (4.3)$$

Indeed, given a function $\varphi$ and a value $\sigma>0$ we may find $m_0\in {\mathbb N}$ such that $m>m_0$ implies $|\varphi(x)-\varphi(y)|\le \sigma/2$ for all $x,y$ such that $\rho(x,y)\le \delta_m$. Take $L=\max_X |\varphi|$ and select $j_0\in {\mathbb N}$ so big that $2L/s_j<\sigma/2$ for any $j>j_0$. Then the absolute value of the expression in the second line of Eq. $(4.3)$ does not exceed $\sigma$. Since $\sigma$ can be taken arbitrarily small, Eq.\, $(4.3)$ is satisfied.

Take a point $\bar x\in \mathrm{supp}\, \mu$. By definition, $\mu_*(B)\neq 0$, where $B=B_{\varepsilon/2}(\bar x)$ is an $\varepsilon/2$~-- ball, centered at $\bar x$. The set $\mathrm{supp}\, \mu_*$ is closed and invariant. By [40, Theorem 1.2.7], it contains a minimal subset. Hence we may assume that $\bar x \in \mathrm{M}(X,T).$ 

Since $\mu_*(B)>0$, there exists an $m_0>0$ such that $J_m(\chi_B)=\mu_m(B)>0$ for all $m>m_0$. Here $\chi_B$ is the characteristic function for the set $B$. By definition of $J_m$ we see that the corresponding set $S_m$ has a positive density in ${\mathbb Z}^+$. Lemma 4.1 is proved. $\square$

Now we suppose that the statement of Theorem 3.1 is wrong. Then there exist a constant $\varepsilon>0$, a positive sequence
$\delta_m\to 0$ and a sequence $p_k^m$ of $\delta_m$~-- pseudotrajectories that cannot be $\varepsilon$~-- shadowed in the sense of $(2.1)$.

Take the point $\bar x$ and the ball $B$ that exist by Lemma 4.1. Fix $m>0$ so that $\mu_m(B)>0$. Let the increasing sequence $\mathcal{I}_m=\{i_j\}$ be such that $p_{i_j}^m\in B$ for all $j\in {\mathbb N}$. By definition of $\mu_m$ we may select $\bar x$ so that
$$
N(\Bar x,B)=\limsup\limits_{n \to\infty} \dfrac{\#(\mathcal{I}_m\bigcap [0,n])}{n}>0.
\eqno (4.4)$$

The set $\{k:T^k(\bar x)\in B\}$ is syndetic (see Lemma 2.10). So, there exists $P>0$ such that for any $k \in {\mathbb N}$ there exists an $s\in \{0,\ldots P\}$ such that $T^k(y_s)\in B$. Here $y_j=T^j(\bar x)$, $j=0,\ldots, P$. 

Let $K_s=\{k_n^s\}\subset {\mathcal I}_m$ be sets such that $T^{k_n^s}(y_s)\in B$, $s=0,\ldots, P$. Evidently, 
$$\mathcal{I}_m=\bigcup_{s=0}^P K_s$$ 
and, by virtue of $(4.4)$ at least one of values 
$$a_r=\limsup\limits_{n\to\infty} \dfrac{\#(K_r\bigcap [0,n])}{n}$$ 
is positive. Then we take $y=y_r$.

To finish the proof, it suffices to observe that $p_k^m,T^k(y)\in B$ implies $\rho(p_k^m,T^k(y))< \varepsilon$. This gives a contradiction to our assumptions on pseudotrajectories $p_k^m$. $\square$

\noindent\textbf{Remark 4.2.} For our proof it is crucial that the space $X$ is compact. There is a simple counterexample to the "non-compact"\ version of the theorem: $X={\mathbb R}$, $T=\mbox{\rm id}$, $x_k=\varepsilon k$.

\noindent\textbf{Proposition 4.3.} \emph{Let $X$ be an arcwise connected compact infinite metric space (e.g. a closed Riemannian manifold), $T:X\to X$ be an invertible equicontinuous map. Then there exists an $\varepsilon_0>0$ such that for any $d>0$ there exists a double -- side $d$ -- pseudotrajectory $x_k$ where none of its two-side subsequences $x_{k_n}$, $k_n \to \pm\infty$ as $n\to \pm\infty$ could be $\varepsilon_0$ shadowed by the subsequence $y_{n_k}$ of a trajectory $\{y_k=T^k (y_0): k\in {\mathbb Z}\}$.}

\noindent\textbf{Proof.} Fix an arc $\Gamma \subset X$ that is not a loop. Consider the parametrizing map $\gamma: [0,1]\to \Gamma$, let $y=\gamma(0)$, $z=\gamma(1)$. Fix $\sigma>0$ so small that $\rho(y,z)>2\sigma$.  

Take $\varepsilon>0$ so that $\rho(x,y)<\varepsilon$ implies $\rho (T^n(x),T^n(y))<\sigma$ for all $n\in {\mathbb Z}$. Particularly, $\varepsilon\le \sigma$.  

Fix a $\delta>0$. Take $\kappa>0$ so that $\rho(x,y)<\kappa$ implies $\rho (T^n(x),T^n(y))<\delta$, $n\in {\mathbb Z}$.

Take $N\in {\mathbb N}$ and a finite sequence $x_k$, $k=0,\ldots, N$ such that $x_0=y$, $x_N=z$ and $\rho (x_{k-1},x_k)\le \kappa$ for all $k=1,\ldots, N$. 

Now we define a sequence $\{p_k\}$ by formulae:
$$p_k= \left[
\begin{array}{l}
T^k (y) \qquad \mbox{if} \quad k\le 0, \\
T^k(x_k) \qquad \mbox{if} \quad 0<k<N,\\
T^k(z) \qquad \mbox{if} \quad k\ge N.
\end{array}
\right.
$$
Observe that $\rho(T(p_k),p_{k+1})=\rho(T^{k+1}(p_k),T^{k+1}(p_{k+1}))\le \delta$ for all $k=0,\ldots,N-1$. Hence $\{p_k\}$ is a $\delta$ -- pseudotrajectory.

If there existed a trajectory $\{q_k=T^k(q_0)\}$ such that $\rho(p_k,q_k)\le \varepsilon$ for any $k\in {\mathbb Z}$, we would have 
$$
\rho(y,z)\le \rho(y, q_0)+\rho(q_0,z)= \rho(y, q_0)+\rho(T^{-N}(q_N),T^{-N}(p_N)\le \sigma+\varepsilon\le 2 \varepsilon.$$
This contradiction finishes the proof. $\square$

As an example, one can consider identical mapping or a rotation of the circle.

\noindent\textbf{Example 4.4.} We give an example of a homeomorphism that does not belong to the class $\mathrm W$. Take the unit circle endowed with the angular coordinate $\varphi$ with the flow defined by ODE $\dot \varphi=\sin^2 \varphi$. (Fig.\, 4). Let $T$ be a discretization of the considered flow. Map $T$ has exactly two fixed points: the west end of the circle $O_w=\{\varphi=\pi\}$ and the east one $O_e=\{\varphi=0\}$. Both these fixed points are semistable. Trajectories of $T$ that do not coincide with one of those points, entirely appertain to the "northern"\ or to the "southern"\ semicircle. In spite of this, pseudotrajectories can "jump"\ through those semistable fixed points and, consequently, rotate infinitely many times around the circle. This proves that $T\notin \mathrm W$. The same example illustrates that $\limsup$ cannot be replaced by $\liminf$ in $(3.1)$. Indeed, for the considered system, pseudotrajectories may stay arbitrarily long in a neighborhood of one fixed point and then leave for another one. So, we can spend $10$ steps in a neighborhood of $O_w$, then (after a fixed number of steps, necessary to proceed from $O_w$ to $O_e$), we wait $10^{10}$ steps in $O_e$, then we go to $O_w$ and spend there $10^{10^{10}}$ steps and so on. In this case, all corresponding lower limits are zero, whatever we select as a shadowing trajectory. 

\begin{figure}[ht!]
\begin{center}
\includegraphics*[width=2in]{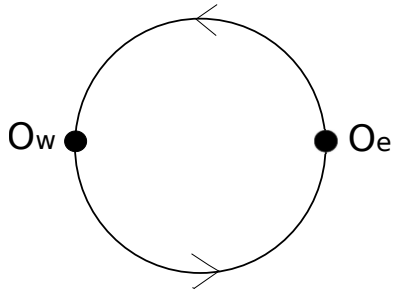}
\end{center}
\caption{No multishadowing for a diffeomorphism of a circle.}
\end{figure}

\noindent\textbf{Corollary 4.5} (to Theorem 3.1). \emph{ For every $\varepsilon>0$ there exists a $\delta>0$ such that for any $\delta$ -- pseudotrajectory $\Xi=\{x_k\}$ of the map $T$ the set $\{k\in {\mathbb Z}: x_k\in U_{\varepsilon} (\overline{\mathrm{M}(X,T))}\}$
is syndetic. Here $U_\varepsilon$ stands for $\varepsilon$ -- neighborhood in the topology of $X$.}

Observe that this statement is very close to one proved by Pilyugin and Sakai [31,32]. The difference is that we take the set $\overline{\mathrm{M}(X,T)}$ instead of $\Omega(X,T)$.

\section{Proof of Lemma 3.5}

Recall that starting from here we always deal with homeomorphisms of compact metric spaces.
Let us prove first of all, that all minimal systems have almost invariant networks.

\noindent{\bf Lemma 5.1.} {\it Any minimal dynamical system $(X,T)$ belongs to the class $Q$.}

\noindent\textbf{Proof}. Take a point $x\in X$. Due to minimality of $(X,T)$ we have $\overline{O(x)}=X$. Fix $\varepsilon>0$, cover $X$ by a finite number of $\varepsilon/2$ -- balls $U_1$, \ldots, $U_K$ and take $n=\max_j n(N(x,U_j))$, see  Eq. $(2.2)$ and Lemma 2.10. Here $n(N(x,U_j))$ is the constant that exists by definition of syndetic sets (Definition 2.9) or, in other words, the maximal possible length of a chain
$$\{T^k(x),T^{k+1}(x),\ldots,T^{m-1}(x),T^m(x)\}\subset X\setminus U_j.$$ Then points $x$, $T(x)$, \ldots, $T^n(x)$ form an almost invariant $\varepsilon$ -- network. $\square$

Now we return to the proof of Lemma 3.5. Consider a finite set $\{b_i\}$ that is an $\varepsilon/2$-network in $X$. By the Bronstein condition, for any $i$ we can find  a minimal set $A_i$ such that $B(b_i,\varepsilon/2) \cap A_i \neq \emptyset$. For every $i$ we select a finite $\varepsilon/2$-network $N_i\subset A_i$ and take $N=\cup_i N_i$. Obviously $N$ is an almost invariant $\varepsilon$-network in $X$. $\square$

\section{Proof of Lemma 3.6} Take an invariant compact subset $K\subset X$ that satisfies conditions of the lemma. Take a point $x_0\in K$ and $\varepsilon>0$. We demonstrate that the set $\overline{B_{\varepsilon}}=\overline{B_{\varepsilon} (x_0)}$ contains a minimal point. Let $B_{\varepsilon/2}=B(\varepsilon/2)(x_0)$.

\noindent\textbf{Lemma 6.1.} \emph{There exists a point $\xi\in {B_{\varepsilon}}$ such that the set $N(\xi,  B_{\varepsilon})$ is syndetic.}

\noindent\textbf{Proof.} Let $A=\{a_i:\, i=1,\ldots, n\}$ be a finite subset of $X$ or "vector". We say that $A$ belongs to the class ${\mathcal H}_\varepsilon$ if for any $k\in {\mathbb Z}$ there exists $j=j(k)$ such that $T^k(a_j)\in 
\overline{B_{\varepsilon/2}}$.  

Observe two evident properties of the class ${\mathcal H}_\varepsilon$.
\begin{enumerate}
\item $A\in {\cal H}_\varepsilon$ if and only if $T^k(A)\in {\cal H}_\varepsilon$ for any $k\in {\mathbb Z}$. 
\item Let $A_k=\{a_i^k,\, i=1,\ldots,n\}\in {\cal H}_\varepsilon,$ $k\in {\mathbb N}$ be a "vector" of $X^n$ converging to a "vector" $A_*$. Then $A_*\in {\cal H}_\varepsilon$.
\end{enumerate}

We start with a set $A=\{a_1,\ldots,a_n\}$ such that $K\subset \bigcup_{i=1}^n T^j(B_{\varepsilon/2}(a_i))$ for any $j\in {\mathbb Z}$. Evidently, $A_*\in {\cal H}_\varepsilon$. Since $x_0\in U_{\varepsilon/2} (A)$,  we may assume that 
$a_1\in B_{\varepsilon/2}$.

Suppose that the set $N(a_1,B_{\varepsilon})$ is non-syndetic (otherwise, we set $\xi=a_1$). Then there exists an increasing sequence $q_m\in {\mathbb N}$ such that $T^{q_m+j}(a_1)\notin B_{\varepsilon}$ for all $j=1,\ldots,m$. Without loss of generality we may suppose that the sequence $T^{q_m} (A)$ converges to a "vector"\ $A_*=\{a_j^*, j=1,\ldots,n\}\in X^n$. Still $A_*\in {\cal H}_\varepsilon$.

Observe that $T^m(a_1^*)\notin B_\varepsilon$ for any $m\in {\mathbb Z}$. Then the $n-1$ point set 
$$A_1=\{a_j^*, j=2,\ldots,n\}$$ 
belongs to the class ${\cal H}_\varepsilon$. Similarly, either the set $N(a_2^*,B_{\varepsilon})$ is syndetic or there exists an $n-2$ point set $A_2\in {\cal H}_\varepsilon$. Repeating this procedure, we must stop after $n$ steps at most and thus obtain the desired point $\xi$. $\square$

Fix the obtained point $\xi$. Let $m\in {\mathbb N}$ be such that the set $N(\xi,\overline{B_{\varepsilon}})=\{n_k\}$ is $m$ -- syndetic. Let $\widetilde\omega$ be the set of all limit points for the sequence $T^{n_k}(\xi)$, $\omega$ be the $\omega$ -- limit set for the trajectory $O(\xi)$. 

Let us prove that
$$ \widetilde{\omega}\subset \overline{B_{\varepsilon}}, \qquad \omega=\widetilde{\omega}\bigcup T(\widetilde{\omega}) \bigcup \ldots \bigcup T^{m} (\widetilde \omega).\eqno (6.1) $$
Indeed, $\widetilde{\omega}\subset \overline{B_{\varepsilon}}$ since $T^{n_k}(\xi)\in \overline{B_{\varepsilon}}$ that is true by definition of $N(\xi,\overline{B_{\varepsilon}})$. Now take a point $\chi \in \omega$. There exists a sequence $p_l$ such that $T^{p_l}(\xi)\to \chi$. Since the set
$N(\xi,\overline B_{\varepsilon})=\{n_k\}$ is $m$ -- syndetic, for any $l\in {\mathbb N}$ we can represent $p_l=n_{k_l}+r_l$ where $r_{l}\in \{0,\ldots m\}$ for all $l$. There is $r\in \{0,\ldots m-1\}$ such that $r_l=r$ for infinitely many values of $l$. We can suppose, proceeding to a subsequence, that $r_l=r$ for all $l$. Then $T^{p_l}(\xi)=T^r(T^{n_{k_l}}(\xi))$ converges to a point of the set $T^r(\widetilde\omega)$. So, $\chi\in T^r(\widetilde\omega)$. 

The set $\omega$ is closed and invariant. Then, by [40,Theorem 1.2.7], it contains a minimal point $\zeta$. By $(6.1)$, there is an iteration $T^q(\zeta)$, $q\in {\mathbb Z}$ that is a point of $\overline U$. This $T^q(\zeta)$ is the desired point. $\square$

\section{Proof of Lemma 3.7} 

Inclusion $\mathrm{Q}\subset \mathrm{W}$ is obvious: iterations of an almost invariant $\varepsilon$--networks trace any sequence, not only pseudotrajectories.

Now suppose that for some $\delta>0$ any $\delta$ -- pseudotrajectory of $T$ is $\varepsilon$ -- multishadowed by a finite set of trajectories. Let us prove existence of an almost invariant $2\varepsilon$ -- network in $\mathrm{CR}(X,T)$.

Consider a point $x \in \mathrm{CR}(X,T)$. Let $\{y_i := y_{i \mod k} | i\in {\mathbb Z}\}$ be a periodic $\delta$ -- pseudotrajectory with $y_0=x$. Due to multishadowing there exists $A(x) := \{a_1 \dots a_k\}$ with $k=k(x)$ such that $x=y_{km}\in B_\varepsilon (T^{km}(A(x))$ for all $m \in {\mathbb Z}$. 

Now select $\{x_1,\ldots, x_N\}$ -- a finite $\varepsilon$ -- network for $\mathrm{CR}(X,T)$. Then
$$A=\bigcup_{j=1}^N \bigcup_{i=0}^{k(x_j)-1} T^i (A(x_j))$$
is such that $\mathrm{CR}(X,T)\subset U_\varepsilon (T^m (A))$ for any $m\in {\mathbb Z}$. 

Now we demonstrate that we can select $A\subset \mathrm{CR}(X,T)$. Take an increasing sequence $\{k_l \in {\mathbb N}\}$ so that iterations $T^{k_l}(A)$ of the set $A$ converge pointwise to a set $A_*$. Then sets $T^m(A_*)\subset \omega(X,T)\subset \mathrm{CR}(X,T)$, $m\in {\mathbb Z}$ form $2\varepsilon$ -- networks there, so it suffices to replace $A$ with $A_*$ and $\varepsilon$ with $2\varepsilon$. $\square$

\section{Proof of Lemma 3.8}

Fix a sequence $\varepsilon_m\to 0$. For every $m$, we consider an almost invariant $\varepsilon_m$ -- network
$$A_m=\{p_{m,j}: j=1,\ldots,N_m\}.$$
Let $\mu_m$ be the probability atomic measure such that $\mu(\{p_{m,j}\})=1/N_m$ for all $j=1,\ldots,N_m$.

Let $T_{\#}$ be the pushforward operator on Borel probability measures induced by $T$. Consider the sequence
$$\mu_{m,n}=\dfrac{1}{n}\sum_{i=0}^{n-1} T_{\#}^i \mu_m.$$
There exists an increasing subsequence $n_l$ such that $\mu_{m,n_l}$ converges in the $*$ -- weak topology. The limit (call it $\mu_m^*$) is a Borel invariant measure. Moreover, for any $x\in X$ we have 
$\mu_{m}^*(B_{\varepsilon_m}(x))\ge 1/N_m$. To construct the desired measure $\mu^*$, we can set
$$\mu^*=\sum_{m=1}^\infty \dfrac{1}{2^m} \mu_m^*.$$
Observe that $U_{\varepsilon_m}(\mathrm{supp}\, \mu_m)=X$ and 
$$\mathrm{supp}\, \mu^*\supset \bigcap_{m=1}^\infty \mathrm{supp}\, \mu_m.$$
So, $\mathrm{supp}\, \mu^*=X$. This finishes the proof. $\square$

\section{Proof of Lemma 3.9}

We start with a statement that is trivial corollary of the definition of chain recurrent sets.

\noindent\textbf{Lemma 9.1.} \emph{For any $\sigma>0$ there exists a $\delta>0$ such that for any $\delta$ -- pseudotrajectory $\Xi=\{x_k\}$ the set 
$$P(X,T,\Xi,\sigma)=\{k\in {\mathbb Z}: x_k\notin U_{\sigma} (\mathrm{CR}(X,T))\}$$
is finite.}

\noindent\textbf{Proof.} Assume that there exists a sequence $\delta_n\to 0$ and a sequence
$$P_n=P(X,T,\Xi_n,\sigma)$$ of infinite sets that correspond to $\delta_n$ -- pseudotrajectories $\Xi_n$. Each of pseudotrajectories $\Xi_n$ has an $\omega$~-- limit point $p_n \notin \overline{U_\sigma({\mathrm{CR}})(X,T)}$. Without loss of generality, we assume that $p_n\to p_*$. Then $p_*\in \mathrm{CR}(X,T)$ that contradicts to our assumptions. $\square$

By (3.3) we have  $\mathrm{CR}(X,T)=\overline{{\mathrm M}(X,T)}$. Bronstein condition implies multishadowing on $\overline{{\mathrm M}(X,T)}$ (Lemmas 3.5 and 3.7).  Given an $\varepsilon>0$ we consider $\delta_0>0$ so that any $\delta_0$~-- pseudotrajectory in $\overline{{\mathrm M}(X,T)}$ is $\varepsilon/2$~-- multishadowed. Take a $\sigma\in (0,\min(\varepsilon/2,\delta_0))$ so that any pointwise $\sigma$~-- perturbation of a $\delta_0/2$~-- pseudotrajectory is a $\delta_0$~-- pseudotrajectory. 

Take $\delta<\delta_0/2$ so that this $\delta$ corresponds to $\sigma$ in the sense of Lemma 9.1. By this lemma any $\delta$~-- pseudotrajectory $p_k$ cannot have infinitely many points out of $\sigma$~-- neighborhood of the set
$\overline{{\mathrm M}(X,T)}={\mathrm{CR}}(X,T)$. 

Fix a $\delta$ -- pseudotrajectory $\{p_k\}$ and consider the sequence $p_k'$ defined as follows. We set $p_k'=p_k$ if $p_k\notin U_\sigma((X,T))$.  Otherwise, we take a point $p_k'\in \mathrm{M}(X,T)$ such that  $\rho(p_k,p_k')<\sigma$. 
The sequence $\{p_k'\}$ is a $\delta_0$~-- pseudotrajectory that consists of two infinite parts inside $\mathrm{M}(X,T)$ and a finite number of points. Such pseudotrajectory can be $\varepsilon/2$~-- traced by a finite number of exact trajectories.  Since $\sigma<\varepsilon/2$, the pseudotrajectory $\{p_k\}$ is $\varepsilon$ -- traced by same trajectories. $\square$

\noindent\textbf{Remark 9.2.} Similarly to Corollary 4.5, we may prove that in conditions of Lemma 9.1
$$\sup_{\Xi\in {\cal P}} P(X,T,\Xi,\sigma)<+\infty.$$
Here $\cal P$ is the set of all $\delta$ -- pseudotrajectories of $T$, $X$, $T$ and $\sigma$ are fixed.

\section{Proof of Lemma 3.10}

Let $x\in \mathrm{CR}(X,T)$. Then for any $\delta>0$ there is a periodic $\delta$ -- pseudotrajectory 
$$\ldots, x=x_0, x_1,x_2, \ldots x_n=x, x_{n+1}=x_1,\ldots $$ 
where $n$ depends on $\delta$. This pseudotrajectory is $\varepsilon$ -- shadowed by a finite number of trajectories $\{T^k(y_m)\}$.  

There exists a $l\in \{1,\ldots,n\}$ such that $\rho(T^{kn}(y_l),x)\le \varepsilon$ for infinitely many $k$. Then there exists a point $q\in \omega(y_l)$ such that $\rho(q,x)\le \varepsilon$. Since $\varepsilon$ is arbitrary, we have proved that
$x\in \Omega(X,T)$. $\square$

\section{Proof of Theorem 3.11} 

\noindent\textbf{1.} By definition, $\mathrm{supp} \,\mu$ is a closed invariant subset of $X$. Fix an $\varepsilon>0$ and consider  $\delta\in (0,\varepsilon)$ such that $\mu(B_\varepsilon (x)) > \delta$ for all $x\in\mathrm{supp}\, \mu$. Such $\delta$ exists since the set $\mathrm{supp} \,\mu$ is compact. Let $A_\delta$ be a finite $\delta$-network, almost invariant with respect to $\mu$. 

Let us prove that 
$$
\mathrm{supp}\, \mu \subset U_{2\varepsilon}(T^n (A_\delta))
\eqno (11.1)$$
for all $n\in {\mathbb Z}$. If (11.1) is not satisfied there exists an $n\in {\mathbb Z}$ and an $\varepsilon$~-- ball $B_\varepsilon(x_0)$, $x_0\in {\mathrm{supp}\, \mu}$ such that $U_\varepsilon(T^n(A))$ for all $n\in {\mathbb Z}$.  Then by definition of almost invariant networks, $\mu (B_\varepsilon (x))\le 1 - (1-\delta)=\delta$. This contradicts to the choice of $\delta$. 

By Lemma 3.6, any neighborhood of any point of $\mathrm{supp}\, \mu$ contains a minimal point. So, $\mathrm{supp}\, \mu\subset \overline{{\mathrm{M}(X,T)}}$.

\noindent\textbf{2.} If minimal points are dense in $\mathrm{supp}\, \mu$, almost invariant $\varepsilon$ -- networks exist by Theorem 3.3. Of course, they all are also almost invariant with respect to $\mu$. 
$\square$

\section{Multishadowing is $C^1$ -- generic}

Certainly, multishadowing is $C^0$ -- generic though the "regular"\ shadowing  is. For instance, multishadowing is $C^0$ -- generic for homeomorphisms of a compact manifold [30].

Here we formulate an important corollary that demonstrates a principle difference between multishadowing and classical shadowing which is not $C^1$ generic [7].

\noindent{\textbf{Theorem 12.1.} \emph{Let $X$ be a $C^1$ smooth compact manifold, $\mathrm {Diff}^1(X)$ be the space of $C^1$ diffeomorphisms. Then the set $W\bigcap \mathrm{Diff}^1(X)$ contains a residual subset in $\mathrm{Diff}^1(X)$.}

\noindent\textbf{Proof.} Given a diffeomorphism $T$, let ${\mathrm P}(X,T)$ be the set of all periodic points. Bonatti and Crovisier [6] demonstrated that for a $C^1$ generic diffeomorphism $T$
$$
\overline{{\mathrm P}(X,T})=\mathrm{CR} (X,T).
\eqno (12.1) $$
By Theorem 3.3, Eq. (12.1) implies that $(X,T)\in \mathrm W$.
$\square$

\section{Discussion}

Let us discuss possible theoretical applications of obtained results: Theorems 3.1 and 3.3. 

We start with Theorem 3.1. Its main idea is quite simple: even an incorrectly applied numerical method can give a correct information about the dynamical system.

Fix a homeomorphism $T$ of a compact metric space $X$.  First of all, recall Corollary 4.5. It claims that any pseudotrajectory has a syndetic set of numbers that correspond to points of the pseudotrajectory in a neighborhood of the set of minimal point.

Basing on technique of Theorem 3.1 we prove that for a sufficiently precise peudotrajectory almost all points are near the set of recurrent points.

\noindent\textbf{Corollary 13.1. }\emph{ For any $\varepsilon>0$ there exists a $\delta>0$ such that for any $\delta$ -- pseudotrajectory $p=\{p_k\}$ of the map $T$
$$\liminf\limits_{N\to \infty} \dfrac{\# K_\varepsilon\bigcap [0,N]}{N}>1-\varepsilon. \eqno (13.1)$$
Here $K_\varepsilon=\{k\ge 0: p_k\in U_\varepsilon ({\mathbf R}(X,T))\}$ 
where $U_\varepsilon({\mathbf R}(X,T))$ is the $\varepsilon$ -- neighborhood of all recurrent points in $X$.}

\noindent\textbf{Proof.} Take a sequence $\delta_m\to 0$ and a sequence $\{p_k^m\}$ of $\delta_m$ pseudotrajectories.  We demonstrate that every $\varepsilon>0$ 
$$\lim\limits_{m\to\infty}\limsup\limits_{N\to \infty} \dfrac{\# L_{m,\varepsilon}\bigcap [0,N]}{N}=0. \eqno (13.2)$$
Here $L_{m,\varepsilon}$ is the completion of the corresponding set $K_{\varepsilon}$ i.e.
$$L_{m,\varepsilon}=\{k\ge 0: p_k^m\notin U_\varepsilon\}.$$

Evidently, (13.2) implies (13.1).

Suppose that (13.2) is wrong. Then, without loss of generality, we may select the sequence $\{p_k^m\}$ so that there is $\alpha>0$ and increasing integer subsequences $\{N_k^m: k\in {\mathbb N}\}$ such that
$$\lim\limits_{m\to \infty} \dfrac{\# L_{m,\varepsilon}\bigcap [0,N_k^m]}{N_k^m}\ge \alpha. \eqno (13.3)$$

For any $m\in {\mathbb N}$ we take a sequence $n_k^m\to \infty$ ($\{n_k^m: k\in {\mathbb N}\}\subset \{N_k^m: k\in {\mathbb N}\}$) such that the limit
$$J_m(\varphi):=\dfrac{1}{n_k^m}\sum_{k=0}^{n_k^m-1} \varphi(p_k^m)$$
is well-defined for any $\varphi \in C^0(X)$ (see the proof of Theorem 3.1, namely Eq. (4.1)).
 
By Riesz Representation theorem, every functional $J_m$ corresponds to a probability measure $\mu_m$ (see Eq.\, (4.2)). We may assume that the sequence $\mu_m$ weakly converges to a measure $\mu_*$ that is invariant (see proof of Theorem 3.1, Section 4). Then $\mathrm{supp}\, \mu_*\subset \overline{\mathrm R}(X,T)$. On the other hand, (13.3) implies that $\mu_m(X\setminus U_\varepsilon ({\mathrm R}(X,T)))\ge \alpha$ for all $m\in {\mathbb N}$. Taking a test function $\varphi$ such that $\varphi(x)=0$ for all $x\in {\mathrm R}(X,T)$ and $\varphi(x) =1$ if $x\notin  U_\varepsilon ({\mathrm R}(X,T))$, we obtain
$$\alpha \le \int_X \varphi\, d\, \mu_m \to \int_X \varphi \, d \mu_*=0.$$

This contradiction finishes the proof. 
$\square$ 

The result of Theorem 3.3 provides a link between Shadowing Theory, Topological Dynamics and Ergodic Theory.  In order to illustrate this we provide two corollaries of Theorem 3.3 and Lemma 3.8. As far as authors know, these results are new even for systems with shadowing or for structurally stable diffeomorphisms that are particular cases of systems considered below.

\noindent\textbf{Corollary 13.2. }\emph{For any homeomorphism $T$ of a compact topological space $X$, such that $\overline{{\mathbf M}(X,T)}=\Omega(X,T)$ and any countable set $\Phi=\{\phi_i : i\in I\}$ of continuous functions there exists an invariant set $\Xi$, dense in $\Omega(X,T)$ and such that for any $\varphi\in \Phi$ and any $x\in \Xi$ there exists a limit
$$\lim\limits_{n\to\infty} \dfrac{1}{n} \sum_{k=0}^{n-1} \varphi (T^k(x)). \eqno (13.3)$$}

\noindent\textbf{Proof.}  Without loss of generality, proceeding to the dynamics on the nonwandering set, we may assume that $\Omega(X,T)=X$. Then, by Lemma 3.8, there exists an invariant probability measure $\mu$ such that $\mathrm{supp}\, \mu=X$. By Birkhof Ergodic Theorem, for any $i\in I$ there exists a set $\Xi_i$ such that the limit (13.3) exists for $\varphi=\phi_i$ and for any $x\in \Xi_i$.

To finish the proof, it suffices to set $\Xi=\bigcap_{i\in I}\Xi_i$. Observe that $\mu(\Xi)=1$ and, since $\mathrm{supp}\,\mu =X$, the set $\Xi$ is dense in $X$.
$\square$

\noindent\textbf{Corollary 13.3. }\emph{For any diffeomorphism $T$ of a compact manifold $X$, such that $\overline{{\mathbf M}(X,T)}=\Omega(X,T)$ there exists a set $\Psi$, dense in $X$ and such that for any $x\in \Psi$ there exists limit
$$\lim_{n\to\infty} \dfrac{1}n \log \|DT^n(x)\|.$$}

Of course, this is the greatest Lyapunov exponent of the trajectory of $x$. 

\noindent\textbf{Proof}. We construct the measure $\mu$, the same as in the previous proof. Then the desired statement follows from Kingmann Subadditive Ergodic Theorem. $\square$

Observe that we may select $\Xi=\Psi$ where sets $\Xi$ and $\Psi$ are defined by Corollaries 13.2 and 13.3 respectively. Indeed, $\mu (\Xi\bigcap \Psi)=1$.  

\section{Conclusion}.

First of all, we list principal results of our paper.

We have established a result that is a weaker version of shadowing (Theorem 3.1): any one-side pseudotrajectory can be shadowed by an exact trajectory along an increasing sequence of time instants. We may assume that points of this trajectory are minimal. 

Certainly, Theorem 3.3 is one of central results of our paper. It gives necessary and sufficient condition of multishadowing and, respectively, new necessary conditions to classical shadowing. It was proved by Aoki and Hirade [3,Theorem 3.1.2] that shadowing property on the chain recurrent set ${\mathrm CR} (X, T)$ implies $(3.3)$. Our Theorem 3.3 improves the mentioned result. First, even the multishadowing property on ${\mathrm CR}(X,T)$ implies (3.3) and, moreover, the Bronstein condition. Particularly, for systems of the class $\mathrm W$, we have $\Omega(\Omega(X,T),T)=\Omega(X,T)$. Also, there must be a probability invariant measure supported on all $\Omega(X,T)$.

Equalities $(3.2)$ and $(3.3)$ are well-known in Dynamics, particularly in Shadowing Theory and $\Omega$ -- Stability Theory. In [25], the authors showed that the following are equivalent:
\begin{itemize}
\item[(a)] $T$ belongs to the set of diffeomorphisms having the periodic shadowing property,
\item[(b)] $T$ belongs to the set of diffeomorphisms having the Lipschitz periodic shadowing property, and
\item[(c)] $T$ satisfies both Axiom A and the no-cycle condition.
\end{itemize}

For Axiom A diffeomorphisms multishadowing is equivalent to (3.3). This follows from Theorem 3.3.

Finally, we list some open problems, that are interesting for us in the framework of our research and may be considered as farther development of our results.
\begin{enumerate}
\item Generally speaking, the density $a$ in Eq.\,  (3.1) depends on the parameter $\varepsilon$ and may tend to zero as $\varepsilon$ tends to zero.  For which systems $(X,T)$ we can take $a$ greater than a fixed positive constant for all $\varepsilon$ and all pseudotrajectories?
\item What does periodic multishadowing property imply?
\item Is there any "two-side"\ version of Theorem 3.1?
\item What can we say about topological entropy for diffeomorphisms with multishadowing?
\item For applications, it is interesting to study systems for which
$$\sup_{\varepsilon>0} N(\varepsilon)<+\infty$$
(see Definition 2.18). Let $\widetilde W$ be the class of such systems. Observe that identical map does not belong to $\widetilde W$. Is it true that ${\widetilde W}={\cal S}$ i.e. does this "strong"\ multishadowing imply shadowing. Is strong multishadowing generic?
\end{enumerate}

\noindent\textbf{Acknowledgements.} Danila Cherkashin was supported by the Saint-Petersburg Mathematical Society (Rokhlin Prize) and by Banach Center of Polish Academy of Sciences. Sergey Kryzhevich was supported by Russian Foundation for Basic Researches, grant 15-01-03797-a and by the Fulbright Visiting Scholar Program. Both coauthors acknowledge Saint-Petersburg State University 
for a research grant 6.38.223.2014.

They are grateful to Prof. Mikhail Blank, Prof. Sergey Pilyugin, Prof. Anatoly Vershik, Andrey Alpeev, Ali Barzanouni and Mikhail Basok for their comments, remarks and suggestions.


\begin{thebibliography}{99}
\bibitem{krycher1} E.\, Akin, \emph{The general topology of dynamical systems}, Grad. Stud. Math., vol.\, 1, American Mathematical Society, Providence, RI, 1993. MR 1219737 (94f:58041).
\bibitem{krycher2} D.\, V.\, Anosov, \emph{Geodesic flows on closed Riemannian manifolds of negative curvature}, Trudy Mat. Inst. Steklov. \textbf{90} (1967), 209 (Russian).
\bibitem{krycher3} N.\, Aoki, K.\, Hiraide, \emph {Topological theory of dynamical systems}, vol. 52 of North- Holland Math. Library. North-Holland Publ. Co. Amsterdam, 1994.
\bibitem{krycher4} N.\,Bernardes, U.Darji, \emph{Graph theoretic structure of maps of the Cantor space}, Adv. Math. \textbf{231} (2012), 1655--1680.
\bibitem{krycher5} M.\, L.\, Blank, \emph{Metric properties of $\varepsilon$ -- trajectory of dynamical systems with stochastic behavior}, Ergodic Theory and Dynamical Systems, \textbf{8}, 1988, 365--378.
\bibitem{krycher6} C.\, Bonatti, S.\, Crovisier, \emph{R\'ecurrence et g\'en\'ericit\'e}, Invent. Math. \textbf{158} (2004), 33--104 (French).
\bibitem{krycher7} C.\, Bonatti, L.\,G.\,Diaz, G.\, Turcat, \emph{Pas de shadowing lemma pour des dynamiques partiellement hyperboliques}, C. R. Acad. Sci. Paris S\'{e}r. I Math. \textbf{330} (2000), no. 7, 587--592 (French).
\bibitem{krycher8} R.\, Bowen, \emph{Equilibrium States and the Ergodic Theory of Anosov Diffeomorphisms}, Lecture Notes in Math., \textbf{470}, Springer-Verlag, 1975.
\bibitem{krycher9} I.\,U. Bronstein, \emph{Extensions of Minimal Transformation Groups}. Providence: American Mathematical Society, 1988.
\bibitem{krycher10} Ch.\, Conley, \emph{Isolated invariant sets and the Morse index}, CBMS Regional Conference Series in Mathematics, 38. American Mathematical Society, Providence, R.I., 1978, ISBN 0-8218-1688-8.
\bibitem{krycher11} D.\,A.\, Dastjerdi, M.\, Hosseini, \emph{Shadowing with chain transitivity}, Topol. Appl. \textbf{156} (2009) 2193--2195.
\bibitem{krycher12} D.\, A.\, Dastjerdi, M.\, Hosseini, \emph{Sub-shadowings}, Nonlinear Analysis, \textbf{72} (2010), 3759--3766.
\bibitem{krycher13} A.\, Fakhari, F.\,H.\, Gane,\emph{On shadowing: ordinary and ergodic}, J. Math. Anal. Appl., \textbf{364} (2010), 151--155.
\bibitem{krycher14} W. H. Gottschalk and G. A. Hedlund, Topological dynamics . Bull. Amer. Math. Soc. \textbf{61} (1955), no. 6, 584--588.
\bibitem{krycher15} E.\,Glasner, \emph{Classifying dynamical systems by their recurrence properties}, Methods Nonlinear Anal. \textbf{24} (2004), 21--40.
\bibitem{krycher16} E.\, Glasner, B.\, Weiss, \emph{Sensitive dependence on initial conditions}, Nonlinearity, \textbf{6} (1993), 1067--1075.
\bibitem{krycher17} A.\, Katok, B.\, Hasselblatt, \emph{Introduction to the Modern Theory of Dynamical Systems}, Cambridge University Press, 1997.
\bibitem{krycher18} S.\, Bezuglyi, S.\, Kolyada, \emph{Topics in Dynamics and Ergodic Theory}, Part of London Mathematical Society Lecture Note Series, 2003, ISBN: 9780521533652
\bibitem{krycher19} P.\, Kocielnak, M.\, Mazur, \emph{Chaos and the shadowing property}, Topol. Appl. \textbf{154} (2007), 2553--2557.
\bibitem{krycher20} S.G.Kryzhevich, \emph{Shadowing along subsequences for continuous mappings}, Vestnik St. Petersburg University: Mathematics, \textbf{47:3} (2014), 102--104.
\bibitem{krycher21} K.\,Lee, K.\, Sakai, \emph{Various shadowing properties and their equivalence}, Discrete and Continuos Dynamical Systems, \textbf{13:2} (2005), 533--539.
\bibitem{krycher22} J.\,H.\, Mai, X.\, Ye, \emph{The structure of pointwise recurrent maps having the pseudo-orbit tracing property}, Nagoya Math J., \textbf{166} (2002), 83--92.
\bibitem{krycher23} T.\,K.\,S.\, Mootathu, \emph{Implications of pseudo-orbit tracing property for continuous maps on compacta}, Top. Appl. \textbf{158} (2011), 2232--2239.
\bibitem{krycher24} T.\,K.\,S.\, Moothathu, P.\,Oprocha, \emph{Shadowing, entropy and minimal subsystems}, Montash Math,\textbf{172} (2013), 357--378.
\bibitem{krycher25} A.\, V.\, Osipov, S.\,Yu.\, Pilyugin and S.\, B.\, Tikhomirov, \emph{Periodic shadowing and $\Omega$-stability}, Regular and Chaoitc Dynam. \textbf{15} (2010), 404--417.1
\bibitem{krycher26} K.\,J.\,Palmer, \emph{Shadowing in Dynamical Systems: Theory and Applications}, Springer, 2009.
\bibitem{krycher27} K.\,J.\, Palmer, S.\,Yu.\,Pilyugin, S.\,B.\,Tikhomirov, \emph{Lipschitz shadowing and structural stability of flows}, Journ. Differ. Equat., \textbf{252} (2012), 1723--1747.
\bibitem{krycher28} S.\,Yu.\, Pilyugin, {\emph{Shadowing in Dynamical Systems}, Lect. Notes Math., Vol. 1706, Springer-Verlag, 1999.}
\bibitem{krycher29} S.\, Yu.\, Pilyugin, \emph{The Space of Dynamical Systems with the $C^0$--Topology}, Lecture Notes in Math., vol. 1571, Springer - Verlag (1994).
\bibitem{krycher30} S.\,Yu.\,Pilyugin and O.\,B.\,Plamenevskaya: \emph{Shadowing is generic}, Topol. Appl., \textbf{97:3}, 253--266 (1999).
\bibitem{krycher31} S.Yu.Pilyugin and K.Sakai, \emph{$C^0$ transversality and shadowing properties}, Proc. Steklov Math. Inst., \textbf{256} (2007), 290--305 .
\bibitem{krycher32} S.\,Yu.\, Piljugin, K.\, Sakai, \emph{Transversality and Shadowing Properties}, Tr. Mat. Inst. Steklova, \textbf{256} (2007), 305--319.
\bibitem{krycher33} S.\, Yu.\, Pilyugin, S.\, B. Tikhomirov, \emph{Lipschitz Shadowing implies structural stability}, Nonlinearity 23 (2010) 2509--2515 = arXiv:1010.3688.
\bibitem{krycher34} C.\, Robinson, \emph{Stability theorems and hyperbolicity in dynamical systems}, Rocky Mount. J. Math., \textbf{7} (1977), 425--437.
\bibitem{krycher35} K.\,Sakai, \emph{Shadowing property and transversality condition}, Dynamical Systems and Chaos (World Sci., Singapore). 1995. - V. 1, P. 233--238.
\bibitem{krycher36} K. Sakai, \emph{Pseudo orbit tracing property and strong transversality of diffeomorphisms of closed manifolds}, Osaka J. Math, \textbf{31} (1994), 373--386.
\bibitem{krycher37} K.\,Sakai, \emph{Various shadowing properties for positively expansive maps}, Topol.\,Appl. \textbf{131:1} (2003), 15--31.
\bibitem{krycher38} K.\,Sawada, \emph{Extended $f$ -- orbits are approximated by orbits}, Nagoya Math. J., \textbf{79} (1980), 33--45.
\bibitem{krycher39} S.\,Smale,  \emph{Differentiable Dynamical Systems}, Bull. Amer. Math. Soc. \textbf{73} (1967), 747--817.
\bibitem{krycher40} J.\, de Vries, \emph{Topological Dynamical Systems. An Introduction to the Dynamics of Continuous Mappings}. De Gruyter, 2014, ISBN 3110342413, 9783110342413
\bibitem{krycher41} C.-C.\, Yuan, J.\,A.\, Yorke, \emph{An open set of maps for which every point is absolutely nonshadowable}, Proc. Amer. Math. Soc., \textbf{128} (2000), 909--918.
\end{thebibliography}
\end{document}